\documentclass[final,nomarks]{dmtcs-episciences}

\usepackage[utf8]{inputenc}

\usepackage[round]{natbib}
\usepackage{amssymb,amsfonts,amsmath,amsthm,color}
\usepackage{enumerate}
\usepackage{array}
\usepackage{xspace}
\newtheorem{thm}{Theorem}[section]
\newtheorem{proposition}[thm]{Proposition}
\newtheorem{cor}[thm]{Corollary}
\theoremstyle{remark}
\newtheorem{rem}[thm]{Remark}

\author{Benjamin R. Jones}
\title{Minor-closed classes of binary functions}
\affiliation{
	Faculty of Information Technology, Monash University, Australia}
\begin{document}
	\publicationdata
	{vol. 26:3}
	{2024}
	{8}
	{10.46298/dmtcs.12230}
	{2023-09-05; 2023-09-05; 2024-05-31}
	{2024-09-10}
	
	\maketitle
	\begin{abstract}
		\ \newline
		Binary functions are a generalisation of the cocircuit spaces of binary matroids to arbitrary functions. Every rank function is assigned a binary function, and the deletion and contraction operations of binary functions generalise matroid deletion and contraction.
		
		We give excluded minor characterisations for the classes of binary functions with well defined minors, and those with an associated rank function. Within these classes, we also characterise the classes of binary functions corresponding to  matroids, binary matroids and polymatroids by their excluded minors. This gives a new proof of Tutte's excluded minor characterisation of binary matroids in the more generalised space of binary functions. 
	\end{abstract}
	
	\section{Introduction}
	A \emph{binary function} on a finite \emph{ground set} $E$ is a function $f\colon 2^E\to\mathbb{R}$ with $f(\emptyset)=1$.
	
	Binary functions were introduced by \cite{farr1993} as a generalisation of cocircuit spaces of binary matroids. A rank transform assigns a (not necessarily matroid) rank function to any binary function. When $f$ is the indicator function of a binary vector space, then the rank transform on $f$ gives the rank function of the associated binary matroid. Aspects of the theory of binary matroids have been extended to binary functions, including deletion, contraction, duality and Tutte-Whitney polynomials by \cite{farr1993,farr2004,farr2007,farr2007a}.
	
	Many important results in graph and matroid theory are characterisations of classes closed under deletion and contraction by a set of \emph{excluded minors}. In graph theory, results range from the characterisation  of planar graphs by \cite{wagner1937} to the theorem of \cite{robertson2004} that any minor-closed class of graphs can be characterised by a finite set of excluded minors. 
	
	The first excluded minor characterisations for matroids were due to \cite{tutte1958,tutte1958a,tutte1959}, and included the classes of regular and graphic matroids.  
	A proof of Rota's conjecture, that the number of excluded minors for matroids representable over any given finite field is finite, has been announced by \cite{geelen2014}.
	The excluded minors for representability over finite fields have been found for the fields $GF(2)$ by \cite{tutte1958}, $GF(3)$ by \cite{bixby1979} and \cite{seymour1979} and $GF(4)$ by \cite{geelen2000}.
	Unlike graphs, minor-closed classes of matroids can have an infinite number of excluded minors. One of the first examples, proved by \cite{lazarson1958}, is the class of matroids representable over any field.
	
	In this paper, we establish the first excluded minor characterisations of classes of binary functions. We find excluded minor characterisations of binary functions for which no minor is undefined, and of those which have rank functions. We also find excluded minor characterisation for the class of binary functions corresponding to binary matroids.
	
	We identify the classes of binary functions corresponding, via the rank transform, to matroids and polymatroids, and give excluded minor characterisations of these classes. By classifying matroids and binary matroids in the larger class of binary functions, we find a novel proof of Tutte's characterisation of binary matroids among matroids.

	\section{Operations on binary functions}
	For a binary function $f\colon 2^E\to\mathbb{R}$, we have $f(\emptyset)=1$. Unless otherwise stated, we will assume that $f$ is a binary function with ground set $E$. The \emph{order} of a binary function is the size of its ground set. 
	
	Contraction and deletion operations on binary functions were introduced by \cite{farr1993}.
	
	The \emph{contraction} of a set $X\subseteq E$ from a binary function $f$ gives the binary function $f/X\colon 2^{E\backslash X}\to\mathbb{R}$ given by
	\begin{align}
		\label{def_con} f/X(Y)=  f(Y).
	\end{align}  
	The \emph{deletion} of $X\subseteq E$ from $f$ gives the binary function
	$f\backslash X\colon 2^{E\backslash X}\to\mathbb{R}$ given by
	\begin{align}
		\label{def_del} f\backslash X(Y)=   \frac{\sum_{Z
				\subseteq X} f(Y\cup Z)}{\sum_{Z
				\subseteq X} f(Z)}.
	\end{align} 
	
	Note that $f/\emptyset =f$ and, since $f(\emptyset)=1$, $\sum_{Z\subseteq \emptyset} f(Z)=f(\emptyset)=1$, so $f\backslash \emptyset=f$.
	
	We write $f/e$  instead of $f/\{e\}$ when contracting by a single element $e\in E$, and likewise write $f\backslash{ e}$ for $f\backslash\{e\}$. Two binary functions $f$ and $f'$, with ground sets $E$ and $E'$ respectively, are \emph{isomorphic} if there is a bijection $\phi:E\to E'$ such that, for all $X\subseteq E$, $f(X)=f'(\{\phi(x):x\in X\})$.
	
	We define a \emph{minor} of a binary function $f$ to be a binary function obtained  from $f$ by a sequence of contractions and deletions of ground set elements. The deletion and contraction operations commute (see \cite[Lemma 4]{farr2004}) in the sense that, for any disjoint $X,Y\subseteq E$, $f\backslash X\backslash Y=f\backslash Y\backslash X$, $f/X/Y=f/Y/X$ and $f\backslash X/ Y=f/ Y\backslash X$ . A minor of a binary function $f$ is \emph{proper} if it has order less than the order of $f$. 
	
	When $\sum_{Z\subseteq X}f(Z)=0$ the minor $f\backslash X$ is undefined. We call a binary function $f$ \emph{stable} if, for all disjoint $X,Y\subseteq E$, $f/X\backslash Y$ is defined, \textit{i.e.} for every $X\subseteq E$,  $\sum_{Z\subseteq X}f(Z)\neq0$.
	
	A class of stable binary functions is \emph{minor-closed} if it closed under taking minors.

	An \emph{excluded minor} of a minor-closed class $\mathcal{C}$ of binary functions is a binary function $f\notin \mathcal{C}$, such that all proper minors of $f$ are in $\mathcal{C}$. Equivalently, $f$ is an excluded minor of $\mathcal{C}$ if $f\notin \mathcal{C}$ and, for every $e\in E$,  $f\backslash{e},f/e\in\mathcal{C}$.
	
	There is only one binary function of order zero, namely $f:\{\emptyset\}\to\mathbb{R}$ where $f(\emptyset)=1$. 
	
	When $|E|=1$, we define $f_\alpha\colon 2^{E}\to\mathbb{R}$ to be the unique (up to isomorphism) order-one binary function with $f_\alpha(E)=\alpha$.
	
	\section{Stable binary functions and rank}
	The first step to understanding minor-closed classes of binary functions is to understand stable binary functions. It turns out that this class is minor-closed. 

	\begin{proposition}
		The stable binary functions are a minor-closed class. 
	\end{proposition}
	\begin{proof}
	Let $g=f\backslash X /Y$ be a minor of a stable binary function $f$ with ground set $E$, and $X'\subseteq E\backslash(X\cup Y)$. As $f$ is stable, both $\sum_{Z\subseteq X}f(Z)$ and $\sum_{Z\subseteq (X\cup X')}f(Z)$ are non-zero. We then have
	\[\begin{array}{rc>{\displaystyle}l}
		0  &   \neq  & \frac{\sum\limits_{Z\subseteq X\cup X'}f(Z)}{\sum\limits_{Z\subseteq X}f(Z)}\\
		&=& \frac{\sum\limits_{V\subseteq X'}\sum\limits_{W\subseteq X}f(V\cup W)}{\sum\limits_{Z\subseteq X}f(Z)}\\
		&=& \sum\limits_{V\subseteq X'}f\backslash X(V)\\&=&\sum\limits_{V\subseteq X'}f\backslash X/Y(V)\\
		&=& \sum\limits_{V\subseteq X'}g(V).\\
	\end{array}\]
	Hence for any subset $X'$ of the ground set of $g$ we have $\sum_{V\subseteq X'}g(V)\neq 0$, and so $g$ is stable. Therefore any minor of a stable binary function is stable.
	\end{proof}
	
	Since this class is closed under taking minors, we can characterise it in terms of its excluded minors.
	
	\begin{thm}{\label{stable}}
		A binary function is stable if and only if it has no minor isomorphic to $f_{-1}$.
	\end{thm}
	\begin{proof}
	$(\Leftarrow)$ As $f_{-1}(\emptyset)+f_{-1}(E)=1-1=0$, we have that $f_{-1}\backslash E$ is not defined, so $f_{-1}$ is not stable. Therefore, any binary function that has $f_{-1}$ as a minor is not stable.

	$(\Rightarrow)$ Let $f$ be a binary function with ground set $E$ that is not stable. Then there is some minimum size $X\subseteq E$ such that $\sum_{Z\subseteq X} f(Z)=0$.  Note that $X\neq \emptyset$, as $f(\emptyset)=1$. 
	
	Let $e\in X$, and consider the minor $g=f\backslash(X\backslash \{e\})/ (E\backslash X)$. As $X$ is minimal, $\sum_{Z\subseteq X\backslash \{e\}} f(Z)\neq 0$ and so $g$ is defined on ground set $\{e\}$. So we have 
	\begin{displaymath}
		\begin{array}{rc>{\displaystyle}l}
			g(\emptyset)	&	=	&	\frac{\sum\limits_{Z\subseteq X\backslash \{e\}}f(Z)}{\sum\limits_{Z\subseteq X\backslash \{e\}}f(Z)}=1,\end{array}
	\end{displaymath}
	
	and
	\begin{eqnarray}
		g(\{e\})	&	=	&	\frac{\sum\limits_{Z\subseteq X\backslash \{e\}}f(Z\cup\{e\})}{\sum\limits_{Z\subseteq X\backslash \{e\}}f(Z)}\\
		&=& \frac{\sum\limits_{Z\subseteq X} f(Z)-\sum\limits_{Z\subseteq X\backslash \{e\}}f(Z)}{\sum\limits_{Z\subseteq X\backslash \{e\}}f(Z)}\\
		&=& \frac{0-\sum\limits_{Z\subseteq X\backslash \{e\}}f(Z)}{\sum\limits_{Z\subseteq X\backslash \{e\}}f(Z)}\\
		&=&-1.
	\end{eqnarray}
	
	Therefore $g=f_1$ is a minor of $f$.	
\end{proof}

	The \emph{rank transform} $Q$ takes a binary function $f$ to its \emph{rank function $Qf\colon 2^{E}\to\mathbb{R}$}, defined in \cite{farr1993} by
	
	\begin{equation}\label{rank_equation}
		Qf(X):=\log_2\left(\frac{\sum\limits_{Y\subseteq E} f(Y)}{\sum\limits_{Y\subseteq E\backslash X} f(Y)}\right).  
	\end{equation}  
	
	In the case $X=E$, we have $Qf(E)=\log_2\left(\sum_{Y\subseteq E} f(Y)\right)$. Note that $Qf$ can be undefined for some binary functions, as the fraction in (\ref{rank_equation}) may be undefined or nonpositive. We say a binary function $f$ is \emph{rankable} if $Qf$ exists. We remark from \cite[p.\ 270]{farr1993} that if $f$ is rankable then  $Qf(\emptyset)=0$ and that for any function $\rho\colon 2^E\to\mathbb{R}$ with $\rho(\emptyset)=0$ there is a binary function $f$ with $Qf=\rho$.
	
	We will now determine some properties of rankable binary functions.
	
	\begin{proposition}\label{ranklemma}
		A binary function $f$ on ground set $E$ is rankable if and only if, for each $X\subseteq E$,
		\[  \sum\limits_{Y\subseteq X} f(Y) >0.    \]
	\end{proposition}
	\begin{proof}
		$(\Rightarrow)$ Suppose $f$ is rankable. Then for each $X\subseteq E$ 
		\[\begin{array}{rc>{\displaystyle}l}
			Qf(E)-Qf(E\backslash X)  &   =  & \log_2\left(\frac{\sum\limits_{Y\subseteq E} f(Y)}{f(\emptyset)}\right)-\log_2\left(\frac{\sum\limits_{Y\subseteq E} f(Y)}{\sum\limits_{Y\subseteq X} f(Y)}\right) \\
			& =    & \log_2\left({\sum\limits_{Y\subseteq X} f(Y)}\right)
		\end{array}\]
		exists, and hence $\sum\limits_{Y\subseteq X} f(Y) >0$.

		$(\Leftarrow)$ Suppose for each $X\subseteq E$ that $\sum\limits_{Y\subseteq X} f(Y) > 0$. Then, for any $Z\subset E$, 
		\[  Qf(Z)= \log_2\left(\frac{\sum\limits_{Y\subseteq E} f(Y)}{\sum\limits_{Y\subseteq E\backslash Z} f(Y)}\right)\]
		is defined, so $f$ is rankable.
	\end{proof}

	\begin{proposition}\label{pro:rankableisstable}
		If $f$ is rankable, then it is stable.
	\end{proposition}
	\begin{proof}
		 We prove the contrapositive.
	Suppose $f\colon 2^E\to\mathbb{R}$ is a non-stable binary function. By Theorem \ref{stable}, $f$ has $f_{-1}$ as a minor, and so $E$ can be partitioned as $E=\{e\}\cup X\cup Z$ such that $f\backslash{X}/Z=f_{-1}$ (on ground set $\{e\}$). Hence we have
	
	\[
	\begin{array}{rl}
		-1=f\backslash{X}/Z(\{e\}) &    =\dfrac{\sum_{Y\subseteq X}f(Y\cup\{e\})}{\sum_{Y\subseteq X}f(Y)} \\&\\
		&=\dfrac{\sum_{Y\subseteq X\cup\{e\}}f(Y)}{\sum_{Y\subseteq X}f(Y)}-1,
	\end{array}\]
	and so $\sum_{Y\subseteq X\cup\{e\}}f(Y)=0$. But by Proposition \ref{ranklemma} we have that $f$ is not rankable.
\end{proof}
	\begin{proposition} \label{pro:rank_minor_closed}
		The class of rankable binary functions is minor-closed.
	\end{proposition}
	
	\begin{proof}
	By Proposition \ref{pro:rankableisstable}, every rankable binary function is stable. We only have to show that rankable binary functions are closed under taking minors.
	
	Let $f\colon 2^E\to\mathbb{R}$ be a rankable binary function, and $e\in E$. Consider contraction. For all $X\subseteq E\backslash \{e\}$,
	\[  \sum\limits_{Y\subseteq (E\backslash \{e\})\backslash X} f/e(Y)= \sum\limits_{Y\subseteq E\backslash (X\cup\{e\})} f(Y)>  0,   \]  
	and so, by Proposition \ref{ranklemma}, $f/e$ is rankable . 
	
	Now consider deletion. We first note that $1+f(\{e\})=\sum_{ Y\subseteq E\backslash( E\backslash\{e\})}f(Y)>0$, so $f\backslash e $ exists. We have, for any $X\subseteq E\backslash\{e\}$,
	
	\[  \sum\limits_{Y\subseteq (E\backslash \{e\})\backslash X} f\backslash e(Y)= \sum\limits_{Y\subseteq E\backslash (X\cup\{e\})} \frac{f(Y)+f(Y\cup \{e\})}{1+f(\{e\})}=\frac{1}{1+f(\{e\})}\sum\limits_{Y\subseteq E\backslash X}f(Y)>  0,   \]  
	and so $f\backslash e$ is rankable. Hence, contraction and deletion of a rankable binary function gives a rankable binary function.
\end{proof}

	We can now find the excluded minors of this class. 
	
	\begin{thm} \label{rank_minors} 
		A binary function is rankable if and only if it has no minor isomorphic to $f_{\alpha}$ for any $\alpha\leq -1$. 
	\end{thm}
	
	\begin{proof} $(\Rightarrow)$ Let $\alpha\leq 1$, and suppose that $f$ has a $f_\alpha$ as a minor.  By Proposition \ref{ranklemma} the rank function of $f_{\alpha}$ is undefined, as $\sum_{Y\subseteq E} f(Y)=1+\alpha\leq 0$. As rankable binary functions are minor-closed by Proposition \ref{pro:rank_minor_closed}, $f$ is not rankable.
	
	$(\Leftarrow)$ let $f \colon 2^E\to \mathbb{R}$ be a non-rankable binary function. Suppose for contradiction that $f$ has no $f_{\alpha}$ minor where $\alpha\leq -1$. Furthermore, let us take $f$ with $|E|$ minimal. Therefore, all proper minors of $f$ are rankable.
	
	As $f$ is not rankable, by Proposition \ref{ranklemma} there is some $X\subseteq E$ such that 
	\[  \sum\limits_{Y\subseteq  X} f(Y) \leq 0.    \]
	As $f(\emptyset)=1$, we must have $X\neq \emptyset$, so let  $e\in X$. Consider $f\backslash{e}$, which is rankable. By Proposition \ref{ranklemma}, we have
	
	\begin{equation*} \label{small_boi}
		0<	\sum\limits_{Y\subseteq X\backslash\{e\}}f\backslash{e}(Y)=\sum\limits_{Y\subseteq X\backslash\{e\}} \frac{f(Y)+f(Y\cup\{e\})}{f(\emptyset)+f(\{e\})} =\frac{1}{1+f(\{e\})}\sum\limits_{Y\subseteq X} f(Y).
	\end{equation*} 
	Hence $1+f(\{e\})< 0$, and so $f/(E\backslash\{e\})\cong f_{\alpha}$, for some $\alpha=f(\{e\})< -1$. Therefore, if $f$ is not rankable, it must have $f_{\alpha}$ as a minor for some $\alpha\leq 1$. 
	\end{proof}
	
	Rankable binary functions are of particular interest, as the deletion and contraction operations of binary functions correspond, via the rank transform, to the rank deletion and contraction operations of matroids. From  \cite[Lemma 4.3]{farr1993} we have the following result.
	\begin{proposition} \label{pro:rank_minor}
		Let $f\colon 2^E\to\mathbb{R}$ be a rankable binary function, and let $X\subseteq E$ and $Y\subseteq E\backslash X$. Then
		\[\begin{array}{l}
			Q(f/X)(Y)=Qf(Y\cup X)-Qf(X),\\
			Q(f\backslash X)(Y)= Qf(Y). 
		\end{array}
		\]
	\end{proposition}
	
	Note here that deletion and contraction for these rank functions are the same operations as those used for rank functions of matroids and polymatroids.

	\section{Indicator functions of binary linear spaces}
	
	A binary function is \emph{straight} if it is the indicator function of a binary vector space. A binary function $f$ is straight if and only if $f$ is  $\{0,1\}$-valued and, for all $X,Y\subseteq E$, if $f(X)=f(Y)=1$ then $f(X\Delta Y)=1$. 	A binary function is \emph{crooked} if it is not straight.

	When a binary function $f$ is the indicator function of a binary vector space then $Qf$ is the  rank function of the binary matroid with that cocircuit space, see
	\cite[Cor.\ 2.2]{farr1993} and \cite[\S 6]{farr2007a}.
	Conversely, for every binary matroid $M=(E,\rho)$, there is a straight binary function -- namely, the indicator function of the cocircuit space of $M$ -- such that $\rho=Qf$. This relationship between the rank and the indicator function of the cocircuit space of a binary matroid was the motivation behind the introduction of binary functions and the rank transform $Q$ by \cite{farr1993}.

	Deletion and contraction operations on binary functions were shown to correspond to matroid deletion and contraction in Proposition \ref{pro:rank_minor}.  As binary matroids are closed under taking minors, the straight binary functions are a minor-closed class.
	
	We have the following result on minors of crooked binary functions, which will be useful in finding the excluded minors of straight binary functions.
	
	\begin{thm}\label{thm:01_nonlinear}
		Let $f$ be a $\{0,1\}$-valued binary function. If $f$ is crooked, then it has a crooked proper minor. 
	\end{thm}
	\begin{proof}  As $f$ is crooked, there is some $X,Y\subseteq E$ such that $f(X)=f(Y)=1\text{, and } f(X\Delta Y)=0$. Let us choose such $X$ and $Y$ with minimum union.
	
	Since $f(\emptyset)=1$, $X$ and $Y$ are distinct, so without loss of generality take $e\in X\backslash Y$. We have
	\begin{align}
		\label{eqn:01-1} f\backslash e (X\backslash \{e\})	=
		\frac{1+f(X\backslash \{e\})}{1+f(\{e\})},	\\
		\label{eqn:01-2}    f\backslash e (Y)=
		\frac{1+f(Y\cup\{e\})}{1+f(\{e\})},  \\
		\label{eqn:01-3} f\backslash e ((X\backslash \{e\})\Delta Y)	=
		\frac{f((X\backslash \{e\})\Delta Y)}{1+f(\{e\})}.
	\end{align}
	
	Suppose that $f(X\backslash \{e\})=1$. As $f(Y)=1$, and $|(X\backslash \{e\}) \cup Y| < |X\cup Y|$, then by the minimality of $|X\cup Y|$, we have $f((X\backslash \{e\})\Delta Y)=1$. By (\ref{eqn:01-1}) and (\ref{eqn:01-3}) we have $f\backslash e (X\backslash \{e\})=2/(1+f(\{e\}))$ and  $f\backslash e((X\backslash \{e\})\Delta Y)=1/(1+f(\{e\}))$ respectively. However, there is no $f(\{e\})\in\{0,1\}$ such that $f\backslash e$ is $\{0,1\}$-valued, and so $f\backslash e$ is crooked.
	
	Let us now consider the case where $f(X\backslash \{e\})=0$. Suppose for contradiction that $f\backslash e$ is straight.
	
	As $f\backslash e$ is $\{0,1\}$-valued,  $f(\{e\})=0$ by (\ref{eqn:01-1}), as  $f(X\backslash \{e\})=0$. This in turn implies that $f(Y\cup\{e\})= 0$ by (\ref{eqn:01-2}).
	
	So by (\ref{eqn:01-1}) and (\ref{eqn:01-2}), $f\backslash e(X\backslash \{e\})=f\backslash e(Y)=1$, and so $f\backslash e((X\backslash \{e\} )\Delta Y)=1$ as $f\backslash e$ is straight. Hence by (\ref{eqn:01-3}) $f((X\Delta Y)\backslash \{e\})=1$.
	
	We have $f(\left(X\backslash \{e\}\right)\Delta Y)=f(Y)=1$, and $|\left((X\backslash \{e\})\Delta Y\right) \cup Y|< |X\cup Y|$, and so by the minimality of $X$ and $Y$  
	we have $f(X\backslash \{e\})=1$, which is a contradiction. Therefore $f\backslash e$ is crooked. 
	\end{proof}
	
	Theorem \ref{thm:01_nonlinear} implies that any excluded minor for the class of straight binary functions is not $\{0,1\}$-valued. We now find the excluded minors for the straight binary functions.

	\begin{thm} \label{linear_minors}
		A binary function is straight if and only if it has no minor isomorphic to any of:
		
		\begin{itemize}
			\item $f_\alpha$, where $\alpha\in\mathbb{R}\backslash\{0,1\}$,
			\item $g\colon 2^{\{a,b\}}\to\mathbb{R}$, where
			\[g(X)=	\left\{\begin{array}{rl}	1 &|X|=0,1\\-1&	|X|=2	\end{array}\right. , \text{ and}	\]
			\item $f_{\mathcal{U}_{2,4}}\colon 2^{\{a,b,c,d\}}\to\mathbb{R}$, where
			\[f_{\mathcal{U}_{2,4}}(X)=	\left\{\begin{array}{rl}
				0&|X|=1,2\\
				1&|X|=0,3\\
				-1&	|X|=4.	\end{array}\right. \]
		\end{itemize}
		
	\end{thm}

	\begin{proof}
	There is only one binary function of order 0, and it is straight. Therefore, all crooked binary functions of order 1 will be excluded minors for the class of straight binary functions.
	
	A binary function of order 1, $f_\alpha$, is straight only if and only if $\alpha\in\{0,1\}$, and so $f_\alpha$ is an excluded minor for the class of straight binary functions for all $\alpha\in\mathbb{R}\backslash{\{0,1\}}$.

	We will now let $f$ be an excluded minor for the class of straight binary functions of order at least two. By Proposition \ref{thm:01_nonlinear}, $f$ is not $\{0,1\}$-valued. For any $X\subset E$, $f/(E\backslash X)$ is a proper minor of $f$, and so it is straight. So $f(X)=f/(E\backslash X)(X)\in\{0,1\}$ for all $X\subset E$. As $f$ is not $\{0,1\}$-valued, we must have $f(E)\notin\{0,1\}$.  
	
	Let $e\in E$. Then
	\begin{equation}\label{eqn:linear_minor_E}
		f\backslash e(E\backslash \{e\}) = \frac{f(E\backslash \{e\})+f(E)}{1+f(\{e\})}.
	\end{equation}
	\[	\]
	
	Since $f\backslash e(E\backslash\{e\}), f(E\backslash\{e\}), f({e})\in \{0,1\}$,  and $f(E)\notin \{0,1\}$, we must have $f(E)\in \{-1,2\}$ by (\ref{eqn:linear_minor_E}).
	\vspace{5mm}
	
	\noindent\textbf{Case 1:}  Suppose $f(E)=2$.  By (\ref{eqn:linear_minor_E}), we have $f(\{e\})=1$ and $f(E\backslash \{e\})=0$ for all $e\in E$. Take distinct $d,e\in E$, and consider the straight binary function $f\backslash e$. We have 
	\begin{align}
		\label{eqn:nonlinear_2_1}   f\backslash e (E\backslash \{e\})	=	\frac{f(E\backslash\{e\})+f(E)}{1+f(\{e\})}	=\frac{0+2}{1+1}=  1,\\
		\label{eqn:nonlinear_2_2}     f\backslash e (\{d\})	=	\frac{f(\{d\})+f(\{d,e\})}{1+f(\{e\})}	=\frac{1+f(\{d,e\})}{2},\\
		\label{eqn:nonlinear_2_3}     f\backslash e (E\backslash\{d,e\})	=	\frac{f(E\backslash\{d,e\})+f(E\backslash\{d\})}{1+f(\{e\})}=	\frac{f(E\backslash\{d,e\})+0}{2}.
	\end{align}
	
	As $f\backslash e(\{d\})\in \{0,1\}$, and $f(E)=2$, we must have $|E|\geq 3$ by (\ref{eqn:nonlinear_2_2}).
	We must then have $f(\{d,e\})\in\{0,1\}$, as $\{d,e\}\subset E$, and so $f\backslash e (\{d\})=1$ by (\ref{eqn:nonlinear_2_2}). Likewise, $f\backslash e(E\backslash\{d,e\})=0$ by (\ref{eqn:nonlinear_2_3}), as $f(E\backslash\{d,e\})\in\{0,1\}$.
	
	Using (\ref{eqn:nonlinear_2_1}), we have $f\backslash e (E\backslash \{e\})=f\backslash e (\{d\})=1$, and $f\backslash e ((E\backslash \{e\}) \Delta \{d\})=f\backslash e (E\backslash \{d,e\})=0$. Hence $f\backslash e$ is crooked, which is a contradiction. So $f$ is not an excluded minor for the straight binary functions if $f(E)=2$.

	\vspace{5mm}
	
	\noindent\textbf{Case 2:} Suppose then that $f(E)=-1$. Then for all $e\in E$, $f(E\backslash \{e\})=1$ by (\ref{eqn:linear_minor_E}), as $f\backslash e (E\backslash \{e\})\in\{0,1\}$. As $|E|\geq 2$ we have, for any distinct $d,e\in E$, 
	
	\begin{equation}
		\label{eqn:nonlinear_2_4}
		f\backslash e (E\backslash\{d,e\})
		=	\frac{f(E\backslash\{d,e\})+1}{1+f(\{e\})}.
	\end{equation}	
	
	As $f\backslash{e}$ is straight, $f\backslash e (E\backslash\{d,e\})\in \{0,1\}$, and so $f(E\backslash\{d,e\})=f(\{e\})$, for any distinct $d,e\in E$. Hence if $f(E)=-1$, then $f(X)$ is constant for all sets $X\subseteq E$ of size 1 or $|E|-2$. This constant is either 0 or 1.
	\vspace{5mm}
	
	\noindent\textbf{Subcase 2.1:} Suppose that $f(X)=1$ for all $X\subseteq E$ of size 1 or $|E|-2$.
	If $|E|=2$, then there is a unique binary function $g$ satisfying these conditions, namely
	\[g(X)=	\left\{\begin{array}{rl}	1,&	|X|=0,1,\\ -1, &|X|=2.\end{array}\right. 	\]
	Deletion and contraction of any element gives a straight binary function, and so $g$ is an excluded minor for the class of straight binary functions.
	
	For  $|E|\geq 3$, take distinct $d,e\in E$. We have
	
	\begin{align}
		\label{eqn:-1_10}
		f\backslash e (E\backslash\{d,e\})=\frac{f(E\backslash\{d,e\})+f(E\backslash \{e\})}{1+f(\{e\})}=\frac{1+1}{1+1}	=1,\\
		\label{eqn:-1_1} f\backslash e (\{d\})	=	\frac{f(\{d\})+f(\{d,e\})}{1+f(\{e\})}= \frac{1+f(\{d,e\})}{1+1},\\
		f\backslash e (E\backslash \{e\})	=	\frac{f(E\backslash \{e\})+f(E)}{1+f(\{e\})}	=	\frac{1-1}{1+1}=  0.
	\end{align}
	
	As $|E|\geq 3$, $\{d,e\}\neq E$ and so $f(\{d,e\})\in\{0,1\}$. Therefore, as $f\backslash{ e}$ is $\{0,1\}$-valued,  $f(\{d,e\})=1$ and $f\backslash{e}(\{d\})=1$ by (\ref{eqn:-1_1}). 
	But then $f\backslash{ e}$ is not straight, as $ f\backslash e (E\backslash\{d,e\})= f\backslash e (\{d\})=1$, but $f\backslash e (E\backslash \{e\})=0$. Hence $f$ is not an excluded minor for the class of straight binary functions.
	\vspace{5mm}
	
	\noindent\textbf{Subcase 2.2:} Suppose that $f(X)=0$ for all $X\subseteq E$  of size 1 or $|E|-2$.
	
	Note, there is no binary function of order 2 which satisfies this condition. For  $|E|\geq 3$, take any distinct $d,e\in E$. We have
	\begin{align}
		\label{eqn:-1_0_1} f\backslash e (E\backslash \{e\})	=	\frac{f(E\backslash \{e\})+f(E)}{1+f(\{e\})}	=	f(E\backslash \{e\})-1,\\	
		\label{eqn:-1_0_2} f\backslash e (E\backslash\{d,e\})=\frac{f(E\backslash\{d,e\})+f(E\backslash \{e\})}{1+f(\{e\})}=f(E\backslash \{e\}),\\
		\label{eqn:-1_0_3} f\backslash e (\{d\})	=	\frac{f(\{d\})+f(\{d,e\})}{1+f(\{e\})}= f(\{d,e\}).
	\end{align}

	As $f\backslash{e}$ is $\{0,1\}$-valued and $f(E\backslash \{e\})=1$, we must have $f\backslash{e}(E\backslash \{e\})=0$ for all $e\in E$ by (\ref{eqn:-1_0_1}). Therefore $f\backslash e (E\backslash\{d,e\})=1$ by (\ref{eqn:-1_0_2}). As $f\backslash e$ is straight, we must have $f\backslash e(\{d\})=f(\{d,e\})=0$ by (\ref{eqn:-1_0_3}). Combining these results gives the following evaluations on  $f$
	\begin{equation*}
		f(X)=	\left\{\begin{array}{rl}	0,&	|X|=1,2,|E|-2,\\ 1,&|X|=0,|E|-1,\\ -1, &\phantom{|}X\phantom{|}=E.	\end{array}\right. 
	\end{equation*}
	
	If $|E|=3$, then no such function satisfies these conditions. If $|E|=4$, these conditions define a  unique (up to isomorphism) binary function, which we will call $f_{\mathcal{U}_{2,4}}$. It is straightforward to confirm by (\ref{def_con}) and (\ref{def_del}) that deletion or contraction of any element from  $f_{\mathcal{U}_{2,4}}$ is straight. Hence $f_{\mathcal{U}_{2,4}}$  is an excluded minor for the class of straight binary functions.
	
	Finally, let us suppose that $|E|\geq 5$. Take distinct $a,b,c,d\in E$. 
	
	By (\ref{eqn:-1_0_2})  and  (\ref{eqn:-1_0_3}) we have  
	
	\begin{align}
		\label{eqn:-1_5+_1}     f\backslash d (\{c\})=0,\\
		\label{eqn:-1_5+_2}    
		f\backslash d (E\backslash\{a,d\})=1,\\
		\label{eqn:-1_5+_3} 
		f\backslash d (E\backslash\{b,d\})=1.
	\end{align}
	
	As $f\backslash d$ is straight and $f\backslash d (E\backslash \{a,d\})=f\backslash d (E\backslash \{b,d\})=1$, then $f\backslash d(\{a,b\})=1$, and so
	
	\begin{equation}
		\label{eqn:-1_5+_4}    1= f\backslash d (\{a,b\})= \frac{f(\{a,b\})+f(\{a,b,d\})}{1+f(\{d\})}=f(\{a,b,d\}).
	\end{equation}
	We have $f\backslash c (\{a,b\})=1$ by (\ref{eqn:-1_5+_4}) and $f\backslash c (\{d\})=0$ by (\ref{eqn:-1_5+_1}), so $f\backslash c (\{a,b,d\})=0$ as $f\backslash{c}$ is straight. But,
	\[    0= f\backslash c (\{a,b,d\})=  \frac{f(\{a,b,d\})+f(\{a,b,c,d\})}{f(\emptyset)+f(\{c\})}=\frac{1+f(\{a,b,c,d\})}{1+0},\]
	and so $f(\{a,b,c,d\})=-1$. But as $|E|\geq 5$, we have $f(\{a,b,c,d\})\in\{0,1\}$. Therefore, no other excluded minors exist. 
	\end{proof}
	\section{Matroids and polymatroids}

	For  $k\in\mathbb{N}$, a binary function $f\colon 2^E\to\mathbb{R}$ is \emph{$k$-polymatroidal} if $f$ is rankable and $Qf\colon 2^E\to\mathbb{R}$ satisfies the following:
	
	\begin{itemize}
		\item[\textbf{{R0{)}}}] \label{itm:R0}$ \forall X\subseteq E, ~Qf(X)\in  \mathbb{N}\cup \{0\} , $
		\item[\textbf{{R1{)}}}] \label{itm:R1} $\forall X\subseteq E,~ Qf(X)\leq k|X|,$
		\item[\textbf{{R2{)}}}] \label{itm:R2} $\forall X\subseteq E,\forall e\in E\backslash X,~  Qf(X)\leq Qf(X\cup\{e\}),$
		\item[\textbf{{R3{)}}}] \label{itm:R3} $\forall X\subseteq E,~\forall  a,b\in E\backslash X \text{ and } a\neq b,~  Qf(X)+Qf(X\cup \{a,b\})\leq Qf(X\cup \{a\})+Qf(X\cup \{b\}). $
	\end{itemize}
	
	A function satisfying \textbf{R0}--\textbf{R3} is the rank function of a $k$-polymatroid (see, for example, \cite{oxley1993}), so $k$-polymatroidal binary functions are those with $k$-polymatroid rank functions. Let $\mathcal{M}_k$ denote the class of $k$-polymatroidal binary functions. 
	
	Note that the deletion and contraction operations for polymatroid rank functions are the same as matroid deletion and contraction, and it is well known (cf.\ \cite{bonin2020}) that $k$-polymatroid rank functions are minor-closed. Combined with Proposition \ref{pro:rank_minor}, this gives the following.

	\begin{rem}
		For each $k\in \mathbb{N}$, $\mathcal{M}_k$ is closed under taking minors.
	\end{rem}

	We can now determine the excluded minors of $\mathcal{M}_k$.
	
	\begin{thm}
		A binary function is  $k$-polymatroidal if and only if it has no minor isomorphic to any of:
		
		\begin{itemize}
			\item 
			$f_\alpha$, where           $\alpha\in\mathbb{R}\backslash\{0,1,3,\ldots,2^{k}-1\}$, and
			\item $f_{\alpha,\beta,\gamma}\colon 2^{\{a,b\}}\to\mathbb{R}$, where $0\leq \alpha\leq \beta \leq \gamma\leq k-1$ are integers and
			\[\left.\begin{array}{lcl}	 
				f_{\alpha,\beta,\gamma}(\emptyset)&=&1,\\
				f_{\alpha,\beta,\gamma}(\{a\})&=&2^{\alpha+\gamma-\beta+1}-1, \\
				f_{\alpha,\beta,\gamma}(\{b\})&=&2^{\gamma+1}-1,\\
				f_{\alpha,\beta,\gamma}(\{a,b\})&=&2^{\alpha+\gamma+1}-2^{\gamma+1}-2^{\alpha+\gamma-\beta+1}+1.\\	\end{array}\right.	\]
		\end{itemize}
	\end{thm}
	
	\begin{proof}
		We first consider $k$-polymatroidal binary functions of order 1. As the unique binary function of order 0 is $k$-polymatroidal, any order 1 binary function $f$ is an excluded minor for $\mathcal{M}_k$ if and only if $f\notin\mathcal{M}_k$.  Again, we write $f_\alpha$ as the binary function of order 1 where $f_\alpha(E)=\alpha$, which is unique up to isomorphism. We then have, where defined, that $Qf_\alpha(\emptyset)=0$ and $Qf_\alpha(E)=\log_2(1+\alpha)$.
		
		By Theorem \ref{rank_minors},  $f_\alpha$ is not rankable if $\alpha\leq -1$. 
		
		If	$\alpha$ is not of the form $2^n-1$ for some $n\in\mathbb{N}\cup\{0\}$, then $Qf_\alpha(E)=\log_2(1+\alpha)\notin\mathbb{N}\cup\{0\}$ and so $Qf$ does not satisfy \textbf{R0}.  
		
		If $\alpha>2^k-1$, then $Qf_\alpha$ does not satisfy \textbf{R1} as $Qf_\alpha(E)>\log_2(1+(2^k-1))=k$. 
		
		If $-1<\alpha<0$, then $Qf_\alpha(E)=\log_2(1+\alpha)<0=Qf_\alpha(\emptyset)$, so $Qf_\alpha$ does not satisfy \textbf{R2}. 
		
		As $Qf_\alpha$ has order one it trivially satisfies \textbf{R3}.

		Hence $f_\alpha$ is $k$-polymatroidal  only if $\alpha\in\{0,1,3,\ldots,2^{k}-1\}$. Finally, it is straightforward to check that $f_\alpha$ is $k$-polymatroidal for each $\alpha\in\{0,1,3,\ldots,2^{k}-1\}$. 
		
		Hence the excluded minors of $\mathcal{M}_k$ of order one are $f_\alpha$, where $\alpha\in\mathbb{R}\backslash\{0,1,3,\ldots,2^{k}-1\}$. 
		
		We now consider excluded minors $f$ of order at least two.
		
		By Theorem \ref{rank_minors}, every binary function of order at least two which is not rankable has a proper minor that is not rankable. Hence any excluded minor of $\mathcal{M}_k$ of order at least two is rankable. 
		
		Let $f$ be an excluded minor of $\mathcal{M}_k$ of order at least 2. Then $Qf$ does not satisfy one of \textbf{R0}, \textbf{R1}, \textbf{R2} or \textbf{R3}.

		Suppose $Qf$ does not satisfy \textbf{R0}. Then there is some $X\subseteq E$ such that $Qf(X)\notin\mathbb{N}\cup\{0\}$. 
		
		As $Qf(\emptyset)=0$, we have that $X$ is nonempty. Take $e\in X$, and as $f/ e\in\mathcal{M}_k$, we know $Q(f/ e)(X\backslash \{e\})\in\mathbb{N}\cup\{0\}$. But $Q(f/ e)(X\backslash \{e\})=Qf(X)-Qf(\{e\})$, and so $Qf(\{e\})\notin \mathbb{N}\cup\{0\}$. But, as $|E|\geq 2$, there is some $d\in E\backslash \{e\}$, and $Q(f\backslash d)(\{e\})=Qf(\{e\})\notin\mathbb{N}\cup\{0\}$, and so $f\backslash d\notin \mathcal{M}_k$.

		Hence $Qf$ must satisfy \textbf{R0}.

		Suppose now that $Qf$ does not satisfy \textbf{R1}. Then there exists $X\subseteq E$ such that $Qf(X)>k|X|$. 
		
		Clearly, $X\neq \emptyset$, so let $e\in X$. Consider $f\backslash (E\backslash \{e\})\in\mathcal{M}_k$, which is a proper minor of $f$, and so $Qf(\{e\})=Q(f\backslash (E\backslash \{e\}))(\{e\})\leq k$.
		
		But then $Q(f/e)(X\backslash \{e\})=Qf(X)-Qf(\{e\})> k|X|-k=k|X\backslash \{e\}|$, and so $f/e\notin\mathcal{M}_k$. 
		
		Hence $Qf$ must satisfy \textbf{R1}.
		
		Suppose now that $Qf$ does not satisfy \textbf{R2}.  That is, there exists $X\subset E$ and $e\in E\backslash X$ such that $Qf(X)>Qf(X\cup \{e\})$. Then we have $Q(f/X)(\{e\})=Qf(X\cup\{e\})-Qf(X)<0$, and so $f/X\notin \mathcal{M}_k$. But since $Qf$ satisfies \textbf{R0} and \textbf{R1}, $X \neq \emptyset$. Hence $f/X$ is a proper minor of $f$, so $f$ is not an excluded minor of $\mathcal{M}_k$. 
		
		Hence $Qf$ must satisfy \textbf{R2}.

		Finally, suppose now that $Qf$ does not satisfy \textbf{R3}.  Therefore there exist distinct $a,b\in E$ and $X\subseteq E\backslash\{a,b\}$ such that
		\begin{equation}\label{eqn:poly_excluded_submodular}
			Qf(X)+Qf(X\cup \{a,b\})> Qf(X\cup \{a\})+Qf(X\cup \{b\}).
		\end{equation}
		
		Suppose there is some $c\in X$. Then, as $f/c\in\mathcal{M}_k$ and satisfies \textbf{R3},
		\[ Q(f/c)(X\backslash \{c\})+   Q(f/c)(X\backslash \{c\} \cup\{a,b\})   \leq Q(f/c)((X\backslash \{c\}) \cup\{a\})  + Q(f/c)((X\backslash \{c\}) \cup\{b\}).   \]
		This can be written as 
		\[Qf(X)+   Qf(X\cup\{a,b\})  -2Qf(\{c\}) \leq Qf(X\cup\{a\})  + Qf(X\cup\{b\}) -2Qf(\{c\}),\]
		which is equivalent to \[Qf(X)+   Qf(X\cup\{a,b\})  \leq Qf(X\cup\{a\})  + Qf(X \cup\{b\}).\]
		This contradicts (\ref{eqn:poly_excluded_submodular}), and so we must have $X=\emptyset$, and (\ref{eqn:poly_excluded_submodular}) becomes
		\begin{equation}\label{eqn:poly_excluded_submodular_updated}
			Qf(\emptyset)+Qf(\{a,b\})> Qf(\{a\})+Qf(\{b\}).
		\end{equation}
		
		Suppose $c\in E\backslash{\{a,b\}}$. Then, as $f\backslash c\in\mathcal{M}_k$ and satisfies \textbf{R3},
		\[Q(f\backslash c)(\emptyset)+   Q(f\backslash c)(\{a,b\})  \leq  Q(f\backslash c)(\{a\})  + Q(f\backslash c)(\{b\}),\]
		which is equivalent to
		\[    Qf(\emptyset)+   Qf(\{a,b\})  \leq Qf(\{a\})  + Qf(\{b\}).\]
		This contradicts (\ref{eqn:poly_excluded_submodular_updated}), and so $E\backslash\{a,b\}$ is empty. Therefore, $f$ must have order 2, \textit{i.e.} $E=\{a,b\}$, and \begin{equation}\label{badpoly_submodular}
			Qf(\{a,b\})>Qf(\{a\})+Qf(\{b\}).
		\end{equation}

		We now know that our excluded minor $f$ must have order 2 and its rank function must satisfy \textbf{R0}, \textbf{R1}, \textbf{R2} and (\ref{badpoly_submodular}). As $Qf$ satisfies (\ref{badpoly_submodular}), we know $f\notin \mathcal{M}_k$, so it remains to consider the assumption that all proper minors of $f$ are $k$-polymatroidal.
		
		We know that $Q(f\backslash a)(\{b\})=Qf(\{b\})$ is a nonnegative integer at most $k$ as $Qf$ satisfies \textbf{R0}, \textbf{R1} and \textbf{R2}. As $Q(f\backslash a)$ has order one, it trivially satisfies \textbf{R3}. Hence, $f\backslash a\in\mathcal{M}_k$, and similarly $f\backslash b \in \mathcal{M}_k$.
		
		We know that $Q(f/ a)(\{b\})=Qf(\{a,b\})-Qf(\{a\})$ is a nonnegative integer as $Qf$ satisfies \textbf{R0}, \textbf{R1} and \textbf{R2}. Hence, $Q(f/ a)$ satisfies \textbf{R0} and \textbf{R2}.  As $Q(f/ a)$ has order one, it trivially satisfies \textbf{R3}. We must then add the restriction that $Qf(\{a,b\})-Qf(\{a\})\leq k$ in order that $Q(f/a)$ satisfies \textbf{R1} and so $f/a\in \mathcal{M}_k$. Similarly, we have $Qf(\{a,b\})-Qf(\{b\})\leq k$ as $f/b\in\mathcal{M}_k$.
		
		Without loss of generality, let $Qf(\{a\})\leq Qf(\{b\})$. The necessary and sufficient conditions for $f$ are

		\[f\colon 2^{\{a,b\}}\to\mathbb{R},\text{ where }        \left\{\begin{array}{l}	 Qf(\{a\}),Qf(\{b\}),Qf(\{a,b\})\in \mathbb{N}\cup\{0\},\\Qf(\{a\})\leq Qf(\{b\})\leq k,\\ Qf(\{a\})+Qf(\{b\})<Qf(\{a,b\})\leq k+Qf(\{a\}).	\end{array}\right.	\]
		Setting $Qf(\{a\})=\alpha,Qf(\{b\})=\beta,$ and $Qf(\{a,b\})=\alpha+\gamma+1$, these conditions are equivalent to $\alpha,\beta,\gamma\in\mathbb{N}\cup\{0\}$ and $0\leq \alpha \leq \beta\leq \gamma \leq k-1$. So the excluded minors of $\mathcal{M}_k$ of order two are of the form
		
		\[\left.\begin{array}{lcl}	 
			f_{\alpha,\beta,\gamma}(\emptyset)&=&1,\\
			f_{\alpha,\beta,\gamma}(\{a\})&=&2^{\alpha+\gamma-\beta+1}-1, \\
			f_{\alpha,\beta,\gamma}(\{b\})&=&2^{\gamma+1}-1,\\
			f_{\alpha,\beta,\gamma}(\{a,b\})&=&2^{\alpha+\gamma+1}-2^{\gamma+1}-2^{\alpha+\gamma-\beta+1}+1.\\	\end{array}\right.	\]
		\end{proof}
		
	It is perhaps surprising that the excluded minors of $k$-polymatroidal binary functions have such a relatively simple structure. While there are an infinite number of excluded minors of order one, there are only $\binom{k+2}{3}$ excluded minors of $\mathcal{M}_k$ of order two.

	Polymatroids and $k$-polymatroid rank functions are a generalisation of the matroid rank functions, which are the $1$-polymatroids. The above theorem can be applied to give the following excluded minor characterisation of binary functions with corresponding matroid rank functions.

	\begin{cor}\label{cor:matroidal_minors}
		A binary function is 1-polymatroidal (\emph{matroidal}) if and only if it has no minor isomorphic to any of:
		
		\begin{itemize}
			\item 
			$f_\alpha$, where           $\alpha\in\mathbb{R}\backslash\{0,1\}$, and
			\item $f\colon 2^{\{a,b\}}\to\mathbb{R}$, where
			\[  f(X)=\left\{    \begin{array}{rl}
				1, & |X|=0,1, \\
				-1,& |X|=2.
			\end{array}    \right.    \]
		\end{itemize}
	\end{cor}

	Theorem \ref{linear_minors} gave an excluded minor characterisation of the straight binary functions. As each straight binary function is matroidal, combining this theorem with Corollary \ref{cor:matroidal_minors} give the excluded minor characterisation of straight binary functions within the class of matroidal binary functions.

	\begin{thm} \label{thm:Tutte_reimagined}
		A matroidal binary function is straight if and only if it has no minor isomorphic to $f_{\mathcal{U}_{2,4}}\colon 2^{\{a,b,c,d\}}\to\mathbb{R}$, where
		\[f_{\mathcal{U}_{2,4}}(X)=	\left\{\begin{array}{rl}
			0,&|X|=1,2,\\
			1,&|X|=0,3,\\
			-1,&	|X|=4.	\end{array}\right. \] 
	\end{thm}
	\begin{proof} Follows from Theorem \ref{linear_minors} and Corollary \ref{cor:matroidal_minors}. 
	\end{proof}
	
	The rank functions of the straight binary functions are the binary matroid rank functions, so this result is equivalent to an excluded minor characterisation of binary matroids within the class of matroids. 
	Here, $Qf_{\mathcal{U}_{2,4}}$ is the rank function of the uniform matroid ${\mathcal{U}_{2,4}}$, and so Theorem \ref{thm:Tutte_reimagined} is equivalent to a classical result in structural matroid theory by \cite{tutte1958,tutte1958a}.
	\begin{thm}[Tutte, 1958]
		A matroid is binary if and only if it has no minor isomorphic to ${\mathcal{U}_{2,4}}$.
	\end{thm}

	\section{Further work}
	This paper lays the groundwork for using binary functions as a tool to investigate minor-closed classes of a number of types of structures.
	
	Binary functions can be used to find new or shorter proofs for excluded minor characterisations of classes of graphs, matroids and other types of objects. In this paper we found a new proof of the characterisation of binary matroids among matroids. One candidate for further investigation are the binary polymatroids, which generalise the binary matroids. It is known from \cite{huynh2010} and \cite{oxley2016} that the set of excluded minors of the binary polymatroids among polymatroids is infinite, but they have not been fully described. If the binary functions that represent binary polymatroids have a simple characterisation, then these could be used to investigate the class of binary polymatroids. 
	
	Binary functions admit an infinite number of minor operations by \cite{farr2004,farr2007}, and one could investigate classes of binary functions closed under minor operations other than deletion and contraction. 
	\\~\\
	
	\noindent \textbf{Acknowledgements} I wish to thank the referees for their detailed and constructive comments, and Graham Farr for his continued support and thoughtful insights. 
\nocite{*}
\bibliographystyle{abbrvnat}
\bibliography{binary_function_excluded_minors}

\begin{thebibliography}{18}
\providecommand{\natexlab}[1]{#1}
\providecommand{\url}[1]{\texttt{#1}}
\expandafter\ifx\csname urlstyle\endcsname\relax
  \providecommand{\doi}[1]{doi: #1}\else
  \providecommand{\doi}{doi: \begingroup \urlstyle{rm}\Url}\fi

\bibitem[Bixby(1979)]{bixby1979}
R.~E. Bixby.
\newblock On {R}eid's characterization of the ternary matroids.
\newblock \emph{J. Combin. Theory Ser. B}, 26\penalty0 (2):\penalty0 174--204,
  1979.
\newblock URL \url{https://doi.org/10.1016/0095-8956(79)90056-X}.

\bibitem[Bonin and Chun(2020)]{bonin2020}
J.~E. Bonin and C.~Chun.
\newblock Decomposable polymatroids and connections with graph coloring.
\newblock \emph{European J. Combin.}, 89:\penalty0 103179, 2020.
\newblock URL \url{https://doi.org/10.1016/j.ejc.2020.103179}.

\bibitem[Farr(1993)]{farr1993}
G.~E. Farr.
\newblock A generalization of the {W}hitney rank generating function.
\newblock \emph{Math. Proc. Cambridge Philos. Soc.}, 113\penalty0 (2):\penalty0
  267--280, 1993.
\newblock URL \url{https://doi.org/10.1017/S0305004100075952}.

\bibitem[Farr(2004)]{farr2004}
G.~E. Farr.
\newblock Some results on generalised {W}hitney functions.
\newblock \emph{Adv. in Appl. Math.}, 32\penalty0 (1-2):\penalty0 239--262,
  2004.
\newblock URL \url{https://doi.org/10.1016/S0196-8858(03)00082-4}.
\newblock Special issue on the Tutte polynomial.

\bibitem[Farr(2007{\natexlab{a}})]{farr2007}
G.~E. Farr.
\newblock On the {A}shkin-{T}eller model and {T}utte-{W}hitney functions.
\newblock \emph{Combin. Probab. Comput.}, 16\penalty0 (2):\penalty0 251--260,
  2007.
\newblock URL \url{https://doi.org/10.1017/S0963548306007966}.

\bibitem[Farr(2007{\natexlab{b}})]{farr2007a}
G.~E. Farr.
\newblock Tutte-{W}hitney polynomials: some history and generalizations.
\newblock In \emph{Combinatorics, complexity, and chance}, volume~34 of
  \emph{Oxford Lecture Ser. Math. Appl.}, pages 28--52. Oxford Univ. Press,
  Oxford, 2007.
\newblock URL \url{https://doi.org/10.1093/acprof:oso/9780198571278.003.0003}.

\bibitem[Geelen et~al.(2014)Geelen, Gerards, and Whittle]{geelen2014}
J.~Geelen, B.~Gerards, and G.~Whittle.
\newblock Solving {R}ota's conjecture.
\newblock \emph{Notices Amer. Math. Soc.}, 61\penalty0 (7):\penalty0 736--743,
  2014.
\newblock URL \url{https://doi.org/10.1090/noti1139}.

\bibitem[Geelen et~al.(2000)Geelen, Gerards, and Kapoor]{geelen2000}
J.~F. Geelen, A.~M.~H. Gerards, and A.~Kapoor.
\newblock The excluded minors for {${\rm GF}(4)$}-representable matroids.
\newblock \emph{J. Combin. Theory Ser. B}, 79\penalty0 (2):\penalty0 247--299,
  2000.
\newblock URL \url{https://doi.org/10.1006/jctb.2000.1963}.

\bibitem[Huynh(2010)]{huynh2010}
A.~Huynh.
\newblock Answer to a question ``{R}epresentability of polymatroids over
  ${GF}(2)$''.\newline
\newblock MathOverflow, 2010-09-02.
\newblock URL
  \url{http://mathoverflow.net/questions/37548/representability-of-polymatroids-over-gf2}.

\bibitem[Lazarson(1958)]{lazarson1958}
T.~Lazarson.
\newblock The representation problem for independence functions.
\newblock \emph{J. London Math. Soc.}, 33:\penalty0 21--25, 1958.
\newblock URL \url{https://doi.org/10.1112/jlms/s1-33.1.21}.

\bibitem[Oxley and Whittle(1993)]{oxley1993}
J.~Oxley and G.~Whittle.
\newblock A characterization of {T}utte invariants of {$2$}-polymatroids.
\newblock \emph{J. Combin. Theory Ser. B}, 59\penalty0 (2):\penalty0 210--244,
  1993.
\newblock URL \url{https://doi.org/10.1006/jctb.1993.1067}.

\bibitem[Oxley et~al.(2016)Oxley, Semple, and Whittle]{oxley2016}
J.~Oxley, C.~Semple, and G.~Whittle.
\newblock A wheels-and-whirls theorem for 3-connected 2-polymatroids.
\newblock \emph{SIAM J. Discrete Math.}, 30\penalty0 (1):\penalty0 493--524,
  2016.
\newblock URL \url{https://doi.org/10.1137/140996549}.

\bibitem[Robertson and Seymour(2004)]{robertson2004}
N.~Robertson and P.~D. Seymour.
\newblock Graph minors. {XX}. {W}agner's conjecture.
\newblock \emph{J. Combin. Theory Ser. B}, 92\penalty0 (2):\penalty0 325--357,
  2004.
\newblock URL \url{https://doi.org/10.1016/j.jctb.2004.08.001}.

\bibitem[Seymour(1979)]{seymour1979}
P.~D. Seymour.
\newblock Matroid representation over {${\rm GF}(3)$}.
\newblock \emph{J. Combin. Theory Ser. B}, 26\penalty0 (2):\penalty0 159--173,
  1979.
\newblock URL \url{https://doi.org/10.1016/0095-8956(79)90055-8}.

\bibitem[Tutte(1958{\natexlab{a}})]{tutte1958}
W.~T. Tutte.
\newblock A homotopy theorem for matroids. {I}.
\newblock \emph{Trans. Amer. Math. Soc.}, 88:\penalty0 144--160,
  1958.
\newblock URL \url{https://doi.org/10.2307/1993243}.

\bibitem[Tutte(1958{\natexlab{b}})]{tutte1958a}
W.~T. Tutte.
\newblock A homotopy theorem for matroids. {II}.
\newblock \emph{Trans. Amer. Math. Soc.}, 88:\penalty0 161--174,
  1958.
\newblock URL \url{https://doi.org/10.2307/1993244}.

\bibitem[Tutte(1959)]{tutte1959}
W.~T. Tutte.
\newblock Matroids and graphs.
\newblock \emph{Trans. Amer. Math. Soc.}, 90:\penalty0 527--552, 1959.
\newblock URL \url{https://doi.org/10.2307/1993185}.

\bibitem[Wagner(1937)]{wagner1937}
K.~Wagner.
\newblock \"{U}ber eine {E}igenschaft der ebenen {K}omplexe.
\newblock \emph{Math. Ann.}, 114\penalty0 (1):\penalty0 570--590, 1937.
\newblock URL \url{https://doi.org/10.1007/BF01594196}.

\end{thebibliography}
\label{sec:biblio}	
\end{document}